\documentclass[12pt]{amsart}
\usepackage{graphicx,amssymb}
\usepackage{mathrsfs}
\usepackage{amssymb,amsmath,amsthm,color}
\usepackage{graphicx,mcite}
\usepackage{hyperref}
\usepackage{url}
\usepackage{setspace}

\theoremstyle{definition}

\theoremstyle{remark}

\numberwithin{equation}{section}

\newcommand{\R}{\mathbb R}

\def\TagOnRight

\def\R {\mathbb{R}}

\newcommand{\be}{\begin{equation}}
\newcommand{\ee}{\end{equation}}
\newcommand{\bea}{\begin{eqnarray}}
\newcommand{\eea}{\end{eqnarray}}
\newcommand{\Bea}{\begin{eqnarray*}}
\newcommand{\Eea}{\end{eqnarray*}}
\newcommand{\bt}{\begin{Theorem}}
\newcommand{\et}{\end{Theorem}}

\newcommand{\bpr}{\begin{Proposition}}
\newcommand{\epr}{\end{Proposition}}

\newcommand{\bl}{\begin{Lemma}}
\newcommand{\el}{\end{Lemma}}
\newcommand{\bi}{\begin{itemize}}
\newcommand{\ei}{\end{itemize}}

\newtheorem{Definition}{Definition}[section]
\newtheorem{Theorem}[Definition]{Theorem}
\newtheorem{Lemma}[Definition]{Lemma}
\newtheorem{Proposition}[Definition]{Proposition}

\newtheorem{Remark}[Definition]{Remark}

\begin{document}
\baselineskip16pt

\title[$L^p-L^q$ estimates for maximal operators]{$L^p-L^q$ estimates for maximal operators associated to families of finite type curves in the plane}
\author[]{Ramesh Manna }

\address{School of Mathematical Sciences, National Institute of Science Education and
Research, Bhubaneswar, An OCC of Homi Bhabha National Institute, Jatni 752050, India.}
\address{Homi Bhabha National Institute, Training School Complex, Anushakti Nagar, Mumbai 400094, India}
\email{rameshmanna@niser.ac.in,rameshmanna21@gmail.com}
\date{\today}

\begin{abstract}

We investigate the \(L^p\)-\(L^q\) boundedness of maximal operators associated with averages over families of finite-type curves in the plane. These operators are defined by
\[
\mathbb{M}f(x):=\sup_{1\leq t\leq 2}
\left|\int_{\mathbb{C}}f(x-ty)\rho(y)\,d\sigma(y)\right|,
\]
where \(d\sigma\) denotes the normalized Lebesgue measure on the curve \(\mathbb{C}\), and \(\rho\) is a smooth cutoff function. Let \(\triangle\) denote the closed triangle with vertices
\[
P=\left(\frac25,\frac15\right),\qquad
Q=\left(\frac12,\frac12\right),\qquad
R=(0,0).
\]
We establish that, whenever
\[
\left(\frac1p,\frac1q\right)
\in
\bigl(\triangle\setminus\{P,Q\}\bigr)
\cap
\left\{
\left(\frac1p,\frac1q\right):
\frac1p<\frac2q,\quad
p>\frac{(m+1)q}{q+1}
\right\},
\]
there exists a constant \(B>0\) such that
\[
\|\mathbb{M}f\|_{L^q(\mathbb{R}^2)}
\leq B\|f\|_{L^p(\mathbb{R}^2)}.
\]
As an application, we derive \(L^p\)-\(L^q\) estimates for solutions of higher-order strictly hyperbolic pseudodifferential equations.

\end{abstract}
\maketitle

\section{Introduction} \label{section1}
Let $\mathbb{C}$ be a smooth curve in the plane. Given a function $f,$ continuous and compactly supported, we consider for each $x \in \mathbb{R}^2$ and $t>0,$ the operator 
$$M_tf(x) \, := \, \int_{\mathbb{C}} f(x-ty) \, \rho(y) \, d\sigma(y),$$
where $d\sigma$ denotes the normalized Lebesgue measure over the curves $\mathbb{C}$ and $\rho\in C_c^{\infty}(\mathbb{C})$ is a smooth cut off function in $\mathbb{R}^2.$ 
Consider now the maximal operator given by 
\begin{eqnarray}
\mathcal{M}f(x) \, := \, \sup_{t>0} |M_tf(x)|. \nonumber
\end{eqnarray}

A natural question, we ask is for what range of the exponents $p\in [1,\infty]$ is  the following a priori inequality satisfied:
\begin{eqnarray} \label{1.1}
\|\mathcal{M}f\|_{L^p(\mathbb{R}^2)} \leq \, B \, \|f\|_{L^p(\mathbb{R}^2)}, ~ f \in \mathscr{S}.
\end{eqnarray}
The study of such a maximal operator over dilations of a fixed curves $\mathbb{C} \subset \mathbb{R}^2$ has its beginnings in the circular maximal theorem of Bourgain \cite{B}.
Bourgain showed that when $ \mathbb{C}=\mathcal{S}^1,$ the unit circle, the corresponding maximal operator is bounded on $L^p(\mathbb{R}^2)$ to $L^p(\mathbb{R}^2)$ for $p > 2.$
His proof of the circular maximal theorem  relies more directly on the geometry involved. The relevant geometry information concerns intersections of pairs of thin annuli, (for more details, see \cite{B}). 
Other proof is due to Mockenhaupt, Seeger and Sogge \cite{MSS2} and proof of this result is based on their local smoothing estimates.
These local smoothing estimates, as well as Bourgain's original techniques actually implies that if one modifies the definition so that the supremum is taken over $1 < t < 2,$
then the resulting circular maximal operator  is bounded from $L^p(\mathbb{R}^2)$ to $L^q(\mathbb{R}^2)$ for some $q > p.$
Here, the maximal operator $\mathbb{M}$ is defined by 
\begin{eqnarray} \label{1.2}
\mathbb{M}f(x) \, := \, \sup_{1 \leq t \leq 2} \left|\int_{\mathbb{C}} f(x-ty) \, \rho(y) \, d\sigma(y) \right|.
\end{eqnarray}

Let $\triangle$ be the closed triangle with vertices $P=(\frac{2}{5}, \frac{1}{5}), ~ Q=(\frac{1}{2}, \frac{1}{2}), ~ R=(0, 0).$ In 1997, Schlag \cite{SG} showed that if $\mathbb{C}$ is unit circle, then the inequality \eqref{1.1} holds if $(\frac{1}{p}, \frac{1}{q})$ lies in the interior of $\triangle.$ His result was obtained using the "combinatorial method" of Kolasa and Wolff cf. \cite{KW}. A different proof of this result was later obtained by Schlag and Sogge cf. \cite{SGS}, which was based on some local smoothing estimates. Schlag also showed that except possibly for endpoints, this result is sharp (see, \cite{SG}, \cite{SGS}). Later, in 2002, Sanghyuk Lee \cite{SL} consider the remaining endpoint estimates for the circular maximal operator.

The aim of this paper is to study the boundedness of maximal operator $\mathbb{M}$ associated with finite type curves (see, Definition \ref{D2.1}) from $L^p(\mathbb{R}^2)$ to $L^q(\mathbb{R}^2).$ 
There is a vast literature on maximal and averaging operators over families of lower dimensional curves in the plane. The main issue turns out to be the curvature: roughly speaking, non-flat curves admit non-trivial maximal estimates; whereas flat curves do not.
More generally, curvature condition plays a crucial role in the analysis of $\mathbb{M},$ and also in other problems in harmonic analysis, cf. \cite{EW, S}, \cite{EMS,RM,RM1}. A fundamental and representative positive result in this direction is the $L^p,~p>2$ boundedness of the  Bourgain's circular maximal operator \cite{B}.

In this paper, we shall consider a situation when the curvature is allowed to vanish of finite order on a finite set of isolated points. In this connections, Iosevich  \cite{IOSE} had already shown that: If $\mathcal{C}$ is of finite-type curve which is of finite type $m$ at $a_0.$ Then, the inequality 
$$\|\mathcal{M}f\|_{L^p(\mathbb{R}^2)} \leq \, B_p \, \|f\|_{L^p(\mathbb{R}^2)}$$ holds for $p > m,$ also this result is sharp.
Here, we shall extend this result to $L^p-L^q$ estimates for the corresponding maximal operator $\mathbb{M}$ associated to families of finite type curves in the plane. 

\section{Some Preliminaries and main result}
We shall need the following definition: 
\begin{Definition} \label{D2.1}
Let $\mathbb{C}:I \rightarrow \mathbb{R}^2,$ where $I$ is a compact interval in $\mathbb{R}$ and $\mathbb{C}$ is smooth. We say that $\mathbb{C}$ is of finite type if $\langle (\mathbb{C}(x) - \mathbb{C}(x_0)), \mu \rangle$ does not vanish of infinite order for any $x_0 \in I,$ and any unit vector $\mu \in \mathbb{R}^2.$
\end{Definition}

We shall also need a more precise definition which would specify the order of vanishing at each point. Let $a_0$ denote a point in the compact interval $I.$ Since our curve is finite type, we can always find a smooth function $\gamma,$ such that in a small neighbourhood of $a_0,$ locally we can write $\mathbb{C}(s)=(s, \gamma(s)),$ where $s \in I.$

\begin{Definition}
Let $\mathbb{C}$ be defined as before. Let $\mathbb{C}(s)=(s, \gamma(s)),$ in a small neighbourhood of $a_0.$ We say that $\mathbb{C}$ is of finite type $m$ at $a_0$ if $\gamma^{(k)}(a_0)=0$ for $1 \leq k < m,$ and $\gamma^{(m)}(a_0) \neq 0.$
\end{Definition}

\begin{Remark}
In fact, this condition is equivalent with the curvature of $\mathbb{C}$ not vanishing to infinite order at $a_0.$
\end{Remark}

We shall use $C$ as a constant independent of $j,$ in several times without mention it. Now we state our main result of this paper. 
\begin{Theorem} \label{T2.4}
Let \(\mathbb{C}\) be a finite-type curve in the plane of type \(m\) at a point \(a_0\). For \(1\leq t\leq 2\), define the averaging operator
\[
M_t f(x):=\int_{\mathbb{C}}f(x-ty)\rho(y)\,d\sigma(y),
\]
where \(d\sigma\) denotes the induced Lebesgue measure on \(\mathbb{C}\), and \(\rho\) is a smooth cutoff function supported in a sufficiently small neighbourhood of \(a_0\). The associated maximal operator is defined by
\[
\mathbb{M}f(x):=\sup_{1\leq t\leq 2}|M_t f(x)|.
\]

Let \(\triangle\) denote the closed triangle with vertices
\[
P=\left(\frac25,\frac15\right),\qquad
Q=\left(\frac12,\frac12\right),\qquad
R=(0,0).
\]
Then, for every pair \((p,q)\) satisfying
\[
\left(\frac1p,\frac1q\right)
\in
\bigl(\triangle\setminus\{P,Q\}\bigr)
\cap
\left\{
\left(\frac1p,\frac1q\right):
\frac1p<\frac2q,\quad
p>\frac{(m+1)q}{q+1}
\right\},
\]
there exists a constant \(B>0\) such that
\begin{equation} \label{thm2.1}
\|\mathbb{M}f\|_{L^q(\mathbb{R}^2)}
\leq B\|f\|_{L^p(\mathbb{R}^2)}.
\end{equation}

%
\end{Theorem}

\begin{Remark}
In view of Iosevich's theorem  (see, \cite{IOSE}), the maximal operator $\mathbb{M}$ is also of course bounded when the exponents lie on the half open line connecting $(\frac{1}{m}, \frac{1}{m})$ and $(0,0)$. Following the example of Iosevich, let $f(x)=|x_{2}|^{-\frac{1}{p}}\log (\frac{1}{|x_{2}|})^{-1}\phi(x)$, where $\phi$ is a smooth positive function supported in $B(0,1)$ and $\phi\equiv1$ near the origin (see~page $1$~of \cite{IOSE}). Then $f\in L^{p}$ for $1<p<\infty$. Locally the curve $\mathcal{C}$ is given by $\mathcal{C}(s)=(s,\gamma(s))$ and without loss
of generality assume that $a_0=0$ and $\gamma(0)=-1$. Now observe that for $p\leq m$,   $\mathcal{A}_{t}f(x)=\infty$ for $|x_{2}|=t$ and $x_{2}<0$. This shows that the estimate \eqref{thm2.1} holds only if $p>m.$ On the other hand, unlike for the case
$p = q$, in the case where $p \ne q$ the $L^p$ -norm and the $L^q$ -norm behave in a different way under the scalings.Indeed, let $m\geq 2$. Consider the curve $\mathcal{C}(s)=(s,a+s^{m})$ and rectangles $R_{1}= [-C\delta^{\frac{1}{m}},C\delta^{\frac{1}{m}}]\times[-C\delta, C\delta]$ and $R_{2}= [-\delta^{\frac{1}{m}},\delta^{\frac{1}{m}}]\times [\frac{4}{3},\frac{5}{3}],$ where $\delta>0$ is small and $C(a)>2.$ Take $f=\chi_{R_{1}}$ and observe that  $\widetilde{\mathcal{M}}_{\operatorname{loc}}f(x)\geq \delta^{\frac{1}{m}}\chi_{R_{2}}(x)$. Thus, one can see that a necessary condition for boundedness
of $\mathcal{M}$ is that
\[p\geq (m+1) \frac{q}{q+1}.\]
\end{Remark}

The proof proceeds in three main steps. We first decompose
the averaging operator $M_t$ into dyadic pieces near the
flat point. Next, we apply the method of stationary phase
to analyze the Fourier transform of the curve measure
associated with each dyadic piece. Finally, following
the stretching argument of Iosevich \cite{IOSE},
we rescale the operators to obtain uniform curvature
estimates. Combining this rescaling with the technical
lemma of [\cite{IOSE}, Lemma~1.3], we reduce the problem
to local smoothing estimates for the corresponding
Fourier integral operators; see \cite{SL,SGS}.


In the proof, we use the following interpolation lemma [\cite{SL}, Lemma $2.6$], which is an extension to multilinear operators of a theorem in \cite{BJ}. We denote by $L^{p, r}$ the Lorentz spaces.

\begin{Lemma} \label{L3.7} (An interpolation lemma)
Let $\epsilon_1, \epsilon_2 > 0. $ Suppose that $\{T_j\}$ is a sequence of linear (or sublinear) operators such that for some $1 \leq p_1, ~ p_2 < \infty,$ and $1 \leq q_1, ~q_2 < \infty,$
$$\|T_jf\|_{q_1} \leq M_1 2^{\epsilon_1 j} \, \|f\|_{p_1}, ~\|T_jf\|_{q_2} \leq M_2 2^{-\epsilon_2 j} \, \|f\|_{p_2}. $$
Then $T= \sum T_j$ is bounded from $L^{p, 1}$ to $L^{q, \infty}$ with 
$$\|Tf\|_{L^{q, \infty}} \leq C \, M_1^{\theta} \, M_2^{1- \theta} \, \|f\|_{L^{p,1}}, $$ where $\theta=\frac{\epsilon_2}{(\epsilon_1
+\epsilon_2)},~ \frac{1}{q}=\frac{\theta}{q_1}+\frac{(1-\theta)}{q_2}, ~\frac{1}{p}=\frac{\theta}{p_1}+\frac{(1-\theta)}{p_2}.$
\end{Lemma} 

We shall also use the following stationary phase result (\cite{CDS2}) for curves in $\mathbb{R}^2.$

\begin{Lemma} \label{3.3}
	Let $\mathbb{C}$ be a smooth curve in the plane with non-vanishing Gaussian curvature and $d \mu$ a $C_c^{\infty}$ measure on $\mathbb{C}.$ Then
	$$|\widehat{d\mu}(\xi)| \leq const \, (1+|\xi|)^{-\frac{1}{2}}.$$
	
	Moreover, suppose that $\Gamma \subset \mathbb{R}^2 \setminus \{0\}$ is the cone consisting of all $\xi$ which are normal to some point $x \in \mathbb{C}$ belonging to a fixed relatively compact neighbourhood $N$ of supp $d\mu.$
	Then,
	$$\left(\frac{\partial}{\partial \xi}\right)^{\alpha} \, \widehat{d\mu}(\xi) =O\left(1+|\xi|\right)^{-N} ~ \forall N,$$ if $\xi \notin \Gamma$ and $\widehat{d\mu}(\xi) \, = \, \sum e^{-i<x_j, \xi>} \, a_j(\xi)$ if $\xi \in \Gamma,$
	where the finite sum is taken over all $x_j \in N$ having $\xi$ as the normal and 
	$$\left|\left(\frac{\partial}{\partial \xi}\right)^{\alpha} a_j(\xi)\right| \, \leq \, C_{\alpha} \, (1+|\xi|)^{-\frac{1}{2}-|\alpha|}.$$
\end{Lemma}

\section{Proof of Theorem \ref{T2.4}}
\begin{proof} of Theorem \ref{T2.4}.
We now turn to the details. Since, our curve is finite type m, after perhaps contracting the support of $\rho$ 
 and making a smooth change of
coordinates, we may assume that the flat point corresponds to $s=0$ and
that the curve is parametrized by
\[
\mathbb C(s)=(s,g(s)s^m+c),\qquad g\in C^\infty,
\quad g(0)\ne0.
\]
We use a smooth dyadic partition of unity in the parameter $s$.

Choose a smooth dyadic partition of unity
$\{\phi_k\}_{k\geq 0}$ on
$\operatorname{supp}\rho\setminus\{0\}$,
such that
\[
\sum_{k\geq 0}\phi_k(s)=1
\]
on $\operatorname{supp}\rho\setminus\{0\}$,
where, for $k\geq 1$,
\[
\operatorname{supp}\phi_k
\subseteq
\{s:2^{-k-1}\leq |s|\leq 2^{1-k}\},
\]
and $\phi_0$ is supported away from the origin.
Setting $\rho_k=\rho\phi_k$, we obtain
\[
\rho(s)=\sum_{k\geq 0}\rho_k(s),
\qquad s\neq 0.
\]

%
%
%
Put
\[
M_t^k f(x)=\int_{\mathbb C}f(x-ty)\phi(2^ky_1)\rho(y)\,d\sigma(y),
\qquad
\mathbb M^k f(x)=\sup_{1\le t\le2}|M_t^k f(x)|.
\]
Thus $\mathbb Mf\le\sum_{k\ge0}\mathbb M^k f$.

Define the anisotropic dilation
\[
\delta_k(x_1,x_2)
=(2^{-k}x_1,2^{-mk}x_2),
\qquad U_kf=f\circ\delta_k.
\]
Then the rescaled operator
$\widetilde M_t^k=U_kM_t^kU_k^{-1}$
is associated with the curve
\[
\gamma_k(s)
=\bigl(s,g(2^{-k}s)s^m+2^{mk}c\bigr),
\qquad \frac12\leq |s|\leq1.
\]
Since
$M_t^k=U_k^{-1}\widetilde M_t^kU_k$
and $|\det\delta_k|=2^{-k(m+1)}$,
we obtain
\begin{equation}\label{eq:scaling-norm}
\|M_t^k\|_{L^p\to L^q}
=
2^{k(m+1)(\frac1p-\frac1q)}
\|\widetilde M_t^k\|_{L^p\to L^q}.
\end{equation}
The same identity holds for the associated
maximal operators.


The curves $s\mapsto(s,g(2^{-k}s)s^m)$ have uniformly bounded
$C^N$ seminorms and curvature bounded away from zero on the fixed
support of $\phi$, for all sufficiently large $k$. The finitely many
remaining $k$ are treated separately. Stationary phase and a finite
conic partition give a principal term of the stretched operator of
the form
\begin{equation}\label{eq:stretched-principal}
\widetilde H_j^k f(x,t)
=\frac{2^{-k}}{(2\pi)^2}\int_\Gamma
 e^{i(\langle x,\xi\rangle+tq_k(\xi)+t2^{mk}c\xi_2)}
 \beta(2^{-j}|\xi|)
 \frac{a_k(t\xi)}{(1+t|\xi|)^{1/2}}
 \widehat f(\xi)\,d\xi,
\end{equation}
where $j\ge1$, $q_k$ is smooth and homogeneous of degree one,
$\operatorname{rank}D_\xi^2q_k=1$ on $\Gamma$, and the symbols $a_k$
have uniformly bounded seminorms. The low-frequency terms and
stationary-phase remainders are kept separately. 
Define the corresponding maximal operators by
\[
\widetilde{\mathbb M}_j^k f(x)
=\sup_{1\le t\le2}
|\widetilde H_j^k f(x,t)|,
\qquad
\mathbb M_j^k f(x)
=\sup_{1\le t\le2}|H_j^k f(x,t)|.
\]
Since
$\widetilde H_j^k=U_kH_j^kU_k^{-1}$,
we have
\[
\mathbb M_j^k f
=
U_k^{-1}\bigl(
\widetilde{\mathbb M}_j^k(U_kf)
\bigr).
\]


Define the order-zero Fourier integral operator
\[
P_j^k f(x,t)
=
\frac{1}{(2\pi)^2}
\int_\Gamma
e^{i(\langle x,\xi\rangle+tq_k(\xi)
+t2^{mk}c\xi_2)}
\beta(2^{-j}|\xi|)
a_k(t\xi)\widehat f(\xi)\,d\xi.
\]

We invoke the following local smoothing estimate
for $P_j^k$, with a fixed smooth cutoff in the
time variable:
\begin{equation}\label{eq:local-smoothing-uniform}
\|P_j^k f\|_{L^q(\mathbb R^2\times[1/2,4])}
\le C2^{j(\frac32-\frac6q)}
\|f\|_{L^p(\mathbb R^2)},
\qquad
\frac1p+\frac3q=1,\quad q>\frac{14}{3},
\end{equation}
where the constant $C$ is independent of
$j$ and $k$.

%

%

For completeness, we record the effect of the translation term in
\eqref{eq:stretched-principal}. On $|\xi| \in [2^{j-1}, 2^{j+1}]$,
\[
\left|q_k(\xi)+2^{mk}c\xi_2\right|
\lesssim 2^j(1+2^{mk}|c|).
\]
The one-dimensional multiplicative Sobolev inequality, followed by
H\"older's inequality in $x$, gives
\[
\|\sup_{1\le t\le2}|F(\cdot,t)|\|_q
\lesssim \|F\|_{L^q_{x,t}}
+\|F\|_{L^q_{x,t}}^{1-1/q}
 \|\partial_tF\|_{L^q_{x,t}}^{1/q}.
\]
The translation $x_2\mapsto x_2+t2^{mk}c$ preserves the
space-time $L^q$ norm, but its derivative contributes the factor
$2^{mk}|c|$. Combining this inequality with
\eqref{eq:local-smoothing-uniform} and the order $-1/2$ of the
stationary-phase amplitude yields
\[
\|\widetilde{\mathbb M}_j^k f\|_q
\lesssim 2^{-k(1-\frac mq)}2^{j(1-\frac5q)}\|f\|_p
\]
for fixed $c\ne0$ (with the same upper bound also valid when $c=0$).
By \eqref{eq:scaling-norm}, we have
\begin{equation}\label{3.10}
\quad
\|\mathbb M_j^k f\|_{L^q(\mathbb R^2)}
\le C2^{-k(1+\frac1q-\frac{m+1}{p})}
 2^{j(1-\frac5q)}\|f\|_{L^p(\mathbb R^2)},\quad
\end{equation}
provided $1/p+3/q=1$ and $14/3<q<\infty$.

\begin{Remark}\label{rem:local-smoothing}
The uniformity in \eqref{eq:local-smoothing-uniform} follows from the
uniform curvature and symbol bounds of the rescaled curves. Indeed,
$g(2^{-k}s)=g(0)+R_k(s)$, where
$\|R_k\|_{C^N}\lesssim_N2^{-k}$ on the fixed parameter interval.
Thus $q_k$ is a uniformly controlled smooth perturbation of the
homogeneous phase for $(s,g(0)s^m)$. The large translation
$2^{mk}c$ is separated from $q_k$ before applying local smoothing.
It does not affect space-time norms, although it must be retained
when differentiating in $t$ to estimate the supremum.
\end{Remark}

\medskip
\noindent\textit{Case $m=2$.}
Let $\mathcal A_{j,t}$ be the frequency-localized averaging operator
for the original curve and set
$\mathbb M_jf=\sup_{1\le t\le2}|\mathcal A_{j,t}f|$.
Since the curvature is nonvanishing on the support of $\rho$, the
stationary-phase bounds and Plancherel's theorem imply
\[
\|\mathcal A_{j,t}f\|_2\lesssim2^{-j/2}\|f\|_2,
\qquad
\|\partial_t\mathcal A_{j,t}f\|_2
\lesssim2^{j/2}\|f\|_2,
\]
uniformly in $1\le t\le2$. Applying the multiplicative Sobolev
inequality in $t$ and integrating in $x$ gives
\[
\begin{aligned}
\|\mathbb M_jf\|_2^2
&\lesssim \|\mathcal A_{j,t}f\|_{L^2_{x,t}}^2
+\|\mathcal A_{j,t}f\|_{L^2_{x,t}}
 \|\partial_t\mathcal A_{j,t}f\|_{L^2_{x,t}}\\
&\lesssim\|f\|_2^2.
\end{aligned}
\]
Consequently, for $j\ge1$,
\begin{equation}\label{3.11}
\|\mathbb M_jf\|_{L^2(\mathbb R^2)}
\le C\|f\|_{L^2(\mathbb R^2)}.
\end{equation}
On the other hand, local smoothing and the time-Sobolev inequality
(with no $k$-dependent translation) give
\begin{equation}\label{3.101}
\|\mathbb M_jf\|_{L^q(\mathbb R^2)}
\le C2^{j(1-\frac5q)}\|f\|_{L^p(\mathbb R^2)},
\qquad \frac1p+\frac3q=1,\quad q>\frac{14}{3}.
\end{equation}
Interpolation between \eqref{3.11} and \eqref{3.101} yields
\begin{equation}\label{3.12}
\|\mathbb M_jf\|_{L^q(\mathbb R^2)}
\le C2^{\frac j2(\frac3p-\frac1q-1)}
\|f\|_{L^p(\mathbb R^2)}
\end{equation}
for the corresponding interpolated pairs. The exponent vanishes on
$3/p-1/q=1$, the line through
$P=(2/5,1/5)$ and $Q=(1/2,1/2)$.
For each point in $[P,Q)$, the admissible interpolation region
contains nearby points on opposite sides of this line. Applying
Lemma~\ref{L3.7} to the frequency-localized sublinear operators
$\mathbb M_j$ therefore gives
\begin{equation}\label{3.13}
\|\mathbb Mf\|_{L^{q,\infty}(\mathbb R^2)}
\le C\|f\|_{L^{p,1}(\mathbb R^2)},
\qquad (1/p,1/q)\in[P,Q).
\end{equation}
The low-frequency operator is bounded separately and does not alter
this estimate. Real interpolation with the $L^\infty\to L^\infty$
bound gives the strong-type estimates in the interior below $[P,Q]$.
The diagonal boundary $p=q>2$ follows from the classical maximal
theorem for curves with nonvanishing curvature. Hence
\begin{equation}\label{3.14}
\|\mathbb Mf\|_{L^q(\mathbb R^2)}
\le C\|f\|_{L^p(\mathbb R^2)}
\end{equation}
for the corresponding points in
$\triangle\setminus\{P,Q\}$ with $p>2$.

\medskip

\noindent\textit{Case $m>2$.}
We first use \eqref{3.10} at $P=(2/5,1/5)$, where
$p=5/2$, $q=5$, and the frequency exponent vanishes.
Since the dyadic cutoff is supported where
$|s| \in [2^{-(k-1)},2^{k+1}]$ and the density of the
curve measure is uniformly bounded, we have
$
\int_{\mathcal C}
|\varphi(2^k y_1)\rho(y)|\,d\sigma(y)
\lesssim 2^{-k}.
$
Consequently,
$\|\mathbb M_j^k f\|_\infty
\lesssim 2^{-k}\|f\|_\infty$.
%
%
Interpolating with \eqref{3.10} at $P$ gives
\[
\|\mathbb M_j^k f\|_{L^q(\mathbb R^2)}
\lesssim
2^{-k(1-\frac{2m+1}{q})}
\|f\|_{L^p(\mathbb R^2)},
\qquad \frac1p=\frac2q,\quad q\ge5.
\]
Hence, summation over $k$ yields
\begin{equation}\label{oneend}
\|\mathbb M_j f\|_{L^q(\mathbb R^2)}
\lesssim \|f\|_{L^p(\mathbb R^2)},
\qquad q=2p>2m+1.
\end{equation}

For the second endpoint, we use the frequency-localized
diagonal maximal estimate of Iosevich \cite{IOSE}.
We first verify the uniform geometric conditions
for the rescaled curves. Set
$h_k(s)=g(2^{-k}s)s^m$. Since $g(0)\neq0$,
we may choose the support of $\rho$ sufficiently
small so that
\[
|m(m-1)g(u)+2mu g'(u)+u^2g''(u)|
\ge c_0>0.
\]
Consequently, on the fixed rescaled support
$|s| \in [1/2,2]$,
\[
|h_k''(s)|\ge c_0|s|^{m-2}\gtrsim1,
\qquad
\sup_{k\ge0}\|h_k\|_{C^N}<\infty
\]
for every $N\ge0$. Moreover,
\[
\|h_k-g(0)s^m\|_{C^N}=O_N(2^{-k}).
\]
Thus the rescaled curves have uniformly
nonvanishing curvature, and the associated
stationary-phase expansions have uniform
symbol bounds.

\begin{Lemma}\label{lem:iosevich-frequency}
Let $m>2$ and $p>m$. Then, for sufficiently
small $\eta>0$, the frequency-localized
maximal operators satisfy
\begin{equation}\label{eq:iosevich-localized}
\|\mathbb M_j^k f\|_{L^p(\mathbb R^2)}
\lesssim
2^{-k(1-\frac mp)}
2^{-j\varepsilon_p}\|f\|_{L^p(\mathbb R^2)},
\qquad j\ge1,\quad k\ge0,
\end{equation}
where
\[
\varepsilon_p=
\begin{cases}
\displaystyle \frac14-\frac1{2p}-\eta,
&3<p\le4,\\[2mm]
\displaystyle \frac1{2p}-\eta,
&p\ge4.
\end{cases}
\]
\end{Lemma}

\begin{proof}
In the stretched variables, the stationary-phase
multiplier contains the factor
$2^{-k}(1+t|\xi|)^{-1/2}$.
After separating the translation term, consider
the order-zero oscillatory integral
\[
\mathcal T_j^k f(x,t)
=\psi(t)\int_\Gamma
e^{i(\langle x,\xi\rangle+tq_k(\xi))}
a_k(t,\xi)\beta(2^{-j}|\xi|)
\widehat f(\xi)\,d\xi.
\]
The local smoothing estimates of
\cite{MSS1,MSS2,IOSE} give
\begin{equation}\label{eq:diagonal-smoothing-input}
\|\mathcal T_j^k f\|_{L^p_{x,t}}
\lesssim_{p,\eta}
2^{j(\sigma(p)+\eta)}\|f\|_{L^p},
\end{equation}
where
\[
\sigma(p)=
\begin{cases}
\displaystyle \frac14-\frac1{2p},
&2<p\le4,\\[2mm]
\displaystyle \frac12-\frac{3}{2p},
&p\ge4.
\end{cases}
\]
The same estimates hold for the corresponding
order-zero amplitudes obtained by normalized
time differentiation, uniformly in $k$.

Let $F(x,t)=\psi(t)\widetilde H_j^k f(x,t)$.
By stationary phase and
\eqref{eq:diagonal-smoothing-input},
\[
\|F\|_{L^p_{x,t}}
\lesssim
2^{-k}2^{-j/2}
2^{j(\sigma(p)+\eta)}\|f\|_p.
\]
Since
\[
|q_k(\xi)+2^{mk}c\xi_2|
\lesssim 2^j(1+2^{mk}|c|),
\qquad |\xi|\in[2^{j-1},2^{j+1}],
\]
we also have
\[
\|\partial_tF\|_{L^p_{x,t}}
\lesssim
2^j(1+2^{mk}|c|)
2^{-k}2^{-j/2}
2^{j(\sigma(p)+\eta)}\|f\|_p.
\]
The multiplicative Sobolev inequality yields
\[
\left\|\sup_{1\le t\le2}|F(\cdot,t)|\right\|_p
\lesssim
\|F\|_{L^p_{x,t}}
+\|F\|_{L^p_{x,t}}^{1-1/p}
\|\partial_tF\|_{L^p_{x,t}}^{1/p}.
\]
Consequently,
\[
\begin{aligned}
\|\widetilde{\mathbb M}_j^k f\|_p
&\lesssim
2^{-k}2^{-j/2}
2^{j(\sigma(p)+\eta)}
\bigl[1+2^{j/p}(1+2^{mk}|c|)^{1/p}\bigr]
\|f\|_p\\
&\lesssim
2^{-k(1-\frac mp)}
2^{-j\varepsilon_p}\|f\|_p,
\end{aligned}
\]
where
\[
\varepsilon_p
=\frac12-\frac1p-\sigma(p)-\eta>0.
\]

It remains to control the stationary-phase
remainder. After separating the oscillatory
phase and translation, uniform stationary
phase gives a remainder amplitude
$r_k(t,\xi)\in S^{-3/2}$, with seminorms
independent of $k$. The corresponding
frequency-localized oscillatory kernels
$K_{j,k,t}$ satisfy
\[
\|K_{j,k,t}\|_{L^1(\mathbb R^2)}
\lesssim 2^{-k}2^{-j},
\qquad
\|\partial_tK_{j,k,t}\|_{L^1(\mathbb R^2)}
\lesssim 2^{-k}(1+2^{mk}|c|).
\]
Thus Young's inequality and the multiplicative
Sobolev inequality give
\[
\left\|\sup_{1\le t\le2}
|\mathcal R_j^k f(\cdot,t)|\right\|_{L^p}
\lesssim
2^{-k(1-\frac mp)}
2^{-j(1-\frac1p)}
\|f\|_{L^p}.
\]
Since $1-\frac1p>\varepsilon_p$, the
remainder has stronger frequency decay
than the principal term. The low-frequency
contribution is handled separately by
smooth kernel estimates.
Finally, anisotropic conjugation preserves
the $L^p\to L^p$ operator norm, proving
\eqref{eq:iosevich-localized}.
\end{proof}

\begin{Remark}\label{rem:iosevich-frequency}
The local smoothing estimates above concern
Fourier integral operators satisfying the
cinematic curvature condition
\cite{MSS1,MSS2,gao2023square}. The translation
term does not affect curvature, since
\[
D_\xi^2(q_k(\xi)+2^{mk}c\xi_2)
=D_\xi^2q_k(\xi).
\]
The time-dependent spatial translation is absorbed
in the space--time estimate, while its contribution
is taken into account when differentiating with
respect to $t$.

The sharp local smoothing estimate of
[\cite{gao2023square}, Theorem~1.4] applies when
$p\ge4$, while for $m=3$ and $3<p<4$,
the earlier, nonsharp estimates of
\cite{MSS1,MSS2,IOSE} are sufficient.
Consequently, $\varepsilon_p>0$ for every
$p>m$, provided $\eta>0$ is chosen
sufficiently small.

\end{Remark}

We now apply Lemma~\ref{lem:iosevich-frequency}.
Since $p>m$, summation over $k$ gives
\begin{equation}\label{otherend}
\|\mathbb M_j f\|_{L^p(\mathbb R^2)}
\lesssim
2^{-j\varepsilon_p}\|f\|_{L^p(\mathbb R^2)},
\qquad p>m,\quad \varepsilon_p>0.
\end{equation}

Interpolating \eqref{oneend} and \eqref{otherend},
and varying the endpoint exponents, gives
positive frequency decay in the interior of
the region bounded by
\[
\left\{\left(\frac2q,\frac1q\right):
q>2m+1\right\}
\quad\text{and}\quad
\left\{\left(\frac1p,\frac1p\right):
p>m\right\}.
\]
Equivalently, within $\triangle$, this region
is characterized by
\[
\frac1p<\frac2q,
\qquad
1+\frac1q-\frac{m+1}{p}>0.
\]
The boundary points are treated by the
corresponding endpoint estimates. Therefore,
at every such interior point there exists
$\varepsilon>0$ such that
\begin{equation}\label{3.001}
\|\mathbb M_j f\|_{L^q(\mathbb R^2)}
\lesssim
2^{-\varepsilon j}
\|f\|_{L^p(\mathbb R^2)}.
\end{equation}
Summing over $j\ge1$ and including the
low-frequency and stationary-phase remainder
contributions yields
\[
\|\mathbb M f\|_{L^q(\mathbb R^2)}
\lesssim \|f\|_{L^p(\mathbb R^2)}.
\]
This completes the proof.

\end{proof}

\subsection{Applications to hyperbolic equations}
In this section, we briefly outline an application of the obtained results to the $L^p-L^q$ estimates for solutions to higher order, strictly hyperbolic pseudo- differential operators. The related fixed-time estimate is given in \cite{RMB, SSS}. The analogous local smoothing result for the wave equation (where $\tilde{q}(\xi)=|\xi|$) is given in \cite{MSS2}.
We obtain the following $L^p-L^q$ estimates for solutions to higher order, strictly hyperbolic pseudo- differential operators.

\begin{Theorem}\label{thm5.1}
Let $m\geq 2$, $N=\lfloor m/2\rfloor$, and
$0<c_1<\cdots<c_N$. Consider the homogeneous,
strictly hyperbolic polynomial
\[
P(\xi,\tau)
=
\tau^{m-2N}
\prod_{j=1}^{N}
(\tau^2-c_j^2|\xi|^2).
\]
Let $u$ be the solution of the Cauchy problem
\begin{equation}\label{5.1}
\begin{cases}
P(D_x,D_t)u(x,t)=0,
& (x,t)\in\mathbb R^2\times(0,\infty),\\
u(x,0)=0,\\
\partial_tu(x,0)=g(x),\\
\partial_t^k u(x,0)=0,
& 2\leq k\leq m-1,
\end{cases}
\end{equation}
where $D_x=-i\nabla_x$ and $D_t=-i\partial_t$.

Suppose that $(1/p,1/q)$ belongs to the
admissible exponent region of
Theorem~\ref{thm2.1} for curves of finite
type $2$. Then
\[
\left\|\sup_{1\leq t\leq 2}|u(\cdot,t)|
\right\|_{L^q(\mathbb R^2)}
\leq C\|g\|_{L^p(\mathbb R^2)}.
\]
\end{Theorem}

\begin{proof}
We first assume that $g\in\mathcal S(\mathbb R^2)$.
Taking the Fourier transform in the spatial
variables and solving the resulting ordinary
differential equation, we obtain
\[
\widehat u(\xi,t)
=
\sum_{j=1}^{N}
d_j\frac{\sin(c_jt|\xi|)}
{c_j|\xi|}\widehat g(\xi),
\]
where the constants $d_j$ are uniquely
determined by
\[
\sum_{j=1}^{N}d_jc_j^{2k}
=
\begin{cases}
1,&k=0,\\
0,&1\leq k\leq N-1.
\end{cases}
\]
Indeed, these identities ensure that the
prescribed initial conditions are satisfied.

By the two-dimensional Poisson formula,
\[
u(x,t)=
\sum_{j=1}^{N}
\frac{d_j}{2\pi c_j}
\int_{|y|<c_jt}
\frac{g(x-y)}
{\sqrt{c_j^2t^2-|y|^2}}\,dy.
\]
For $r>0$, define the circular averaging
operator
\[
A_rg(x)=\frac{1}{2\pi}
\int_{S^1}g(x-r\omega)\,d\sigma(\omega).
\]
Passing to polar coordinates and making
the substitution $r=c_jts$, we obtain
\[
u(x,t)
=
t\sum_{j=1}^{N}d_j
\int_0^1
\frac{s}{\sqrt{1-s^2}}
A_{c_jts}g(x)\,ds.
\]
Consequently,
\[
\sup_{1\leq t\leq 2}|u(x,t)|
\leq
2\sum_{j=1}^{N}|d_j|
\int_0^1
\frac{s}{\sqrt{1-s^2}}
\sup_{1\leq t\leq 2}
|A_{c_jts}g(x)|\,ds.
\]

Set
\[
\delta=\frac1p-\frac1q.
\]
Since the unit circle has everywhere
nonvanishing curvature, Theorem~\ref{thm2.1},
together with a finite partition of unity
and scaling, yields
\[
\left\|\sup_{1\leq t\leq 2}
|A_{at}g|\right\|_q
\leq Ca^{-2\delta}\|g\|_p,
\qquad a>0.
\]
Therefore, by Minkowski's inequality,
\[
\begin{aligned}
\left\|\sup_{1\leq t\leq 2}|u(\cdot,t)|
\right\|_q
&\leq
C\sum_{j=1}^{N}|d_j|c_j^{-2\delta}
\int_0^1
\frac{s^{1-2\delta}}{\sqrt{1-s^2}}\,ds
\,\|g\|_p\\
&\lesssim \|g\|_p.
\end{aligned}
\]
The integral is finite since
$\delta<1$ throughout the admissible
region of Theorem~\ref{thm2.1}.
The general case follows by density.
This completes the proof.
\end{proof}

\begin{Remark}
When $m$ is odd, the characteristic polynomial
$P(\xi,\tau)$ has an additional zero root.
Nevertheless, all its characteristic roots
are real and pairwise distinct for
$\xi\neq0$, so strict hyperbolicity is
preserved. Thus Theorem~\ref{thm5.1}
holds for every integer $m\geq2$.
\end{Remark}

\noindent
{\bf Acknowledgements:} 
\noindent
RM acknowledges the partial support provided by the research grants  (ANRF/ARGM/2025/000750/MTR).

\noindent
{\bf AI disclosure.}   AI tools were used to assist in the revision of this manuscript. The author takes full responsibility for its content.


\end{document}